%% file: main.tex
\newif\ifreport
\reporttrue

\ifreport
  \documentclass{article}\def\zibreport{1}

  \usepackage{authblk}
  \makeatletter
  \renewcommand\AB@authnote[1]{\rlap{\textsuperscript{\normalfont#1}}}

  \makeatother
\else
  \documentclass{llncs}\def\zibreport{0}
\fi

\usepackage[T1]{fontenc}
\usepackage[utf8x]{inputenc}
\usepackage[english]{babel}

\usepackage{booktabs}
\usepackage{pdflscape}
\usepackage{tabularx}
\usepackage{longtable}

\usepackage{tikz}
\usetikzlibrary{decorations.pathmorphing}

\usepackage{pgfplots, pgfplotstable}
\pgfplotsset{compat=1.8}

\ifreport
\usepackage[font=small]{caption}
\fi

\newcommand{\thetitle}{A Status Report on Conflict Analysis in Mixed Integer Nonlinear Programming}

\include{Defines}

\ifthenelse{\zibreport = 1}{
  \usepackage{zibtitlepage}

  \ZTPTitle{\thetitle}
  \ZTPAuthor{
    \ZTPHasOrcid{Jakob Witzig}{0000-0003-2698-0767}
    \hfill
    Timo Berthold
    \hfill
    Stefan Heinz
  }
  \ZTPPreprint
  \ZTPNumber{18-57}
  \ZTPMonth{November}
  \ZTPYear{2018}
}{}

\begin{document}

\ifreport
  \author[1]{Jakob~Witzig}
  \author[2]{Timo~Berthold}
  \author[2]{Stefan~Heinz}

  \affil[1]{Zuse Institute Berlin, Takustr.~7, 14195~Berlin, Germany \protect\\ \texttt{witzig@zib.de}\medskip}
  \affil[2]{Fair Isaac Germany GmbH, Takustr.~7, 14195~Berlin, Germany \protect\\ \texttt{\{timoberthold,stefanheinz\}@fico.com}}

  \zibtitlepage
\else
   \author{
     Jakob Witzig\inst{1}
     \and Timo Berthold\inst{2}
     \and Stefan Heinz\inst{2}
   }
   \institute{
     Zuse Institute Berlin, Takustr.~7, 14195~Berlin, Germany\\ \email{witzig@zib.de}
     \and
     Fair Isaac Germany GmbH, Takustr.~7, 14195~Berlin, Germany\\ \email{\{timoberthold,stefanheinz\}@fico.com}
   }
\fi

\title{\thetitle}

\maketitle

\begin{abstract}
  Mixed integer nonlinear programs (MINLPs) are arguably among the hardest optimization problems, with a wide range of applications.
%  and is used if approximating a nonlinear optimization problem
%  by a linear model does not reflects the entire model sufficiently.
  MINLP solvers that are based on linear relaxations and spatial branching work similar as mixed integer programming (MIP)
  solvers in the sense that they
  are based on a branch-and-cut algorithm, enhanced by various heuristics, domain propagation, and presolving techniques.
  However, the analysis of infeasible subproblems, which is an important component of most major MIP solvers,
  has been hardly studied in the context of MINLPs.
  There are two main approaches for infeasibility analysis in MIP solvers: \emph{conflict graph analysis}, which
  originates from artificial intelligence and constraint programming, and \emph{dual ray analysis}.

%  When solving MINLP with a combination of outer approximation and spatial \bandb,
%  most of MINLP solvers use conflict graph analysis as known from MIP.
%  However, several recent computational studies report that applying this well studied
%  technique has almost no affect for solving MINLPs in practice.
  The main contribution of this short paper is twofold.
  Firstly, we present the first computational study regarding the impact of dual ray analysis on convex and nonconvex MINLPs.
%  For our implementation we used the academic constraint integer programming solver \scip.
  In that context, we introduce a modified generation of infeasibility proofs that incorporates linearization cuts
  that are only locally valid.
  Secondly, we describe an extension of conflict analysis that works directly with the nonlinear relaxation of convex MINLPs
  instead of considering a linear relaxation.
  This is work-in-progress, and this short paper is meant to present first
  theoretical considerations without a computational study for that part.
\end{abstract}

\section{Introduction}

In this paper, we consider \emph{mixed integer nonlinear programs} (MINLPs) of the form
% \begin{align}
%     \min c^\T\primalvector & \notag \\
%     s.t.\quad A\primalvector &\geq b \notag \\
%     g_k(\primalvector) &\leq 0 && \forall k \in \nlcons,  \\
%     \lb \leq \primalvector &\leq \ub && \notag \\
%     \primalvector_j &\in \Z && \forall j \in \intvars, \notag
% \end{align}
\begin{align}
    \min \{c^\T\primalvector \;|\; A\primalvector \geq b,\ g_k(\primalvector) \leq 0\; \forall k \in \nlcons,\
    \lb \leq \primalvector \leq \ub,\ \primalvector_j \in \Z\; \forall j \in \intvars\} \label{eq:minlp}
\end{align}
with objective coefficient vector $c \in \R^{n}$, linear constraint matrix $A \in \R^{m \times n}$,
nonlinear constraint functions $g_k \colon \R^n \mapsto \R$, $k \in \nlcons := \{1,\ldots,p\}$, continuously differentiable, and possibly nonconvex,
and variable bounds $\ell,u \in \overline{\R}^{n}$, where $\overline{\R} := \R \cup \{\pm\infty\}$.
Furthermore, let $\allvars = \{1,\ldots,n\}$ be the index set of all variables and
$\intvars \subseteq \allvars$ the set of variables that need to be integral in every feasible solution.
Without loss of generality, we assume the objective function to be linear.
A nonlinear objective function can be transformed into a constraint bounded
by an artificial variable $z$ that needs to be minimized.
We call an MINLP \emph{convex} when all of its constraint functions $g_k$ are convex.
Otherwise, we call the MINLP \emph{nonconvex}.
When omitting the integrality requirements, we obtain the \emph{nonlinear programming} (NLP) relaxation of~\eqref{eq:minlp}
\begin{align}
  \min \{ c^\T \primalvector\,|\,A\primalvector \geq b,\ g_k(\primalvector) &\leq 0\, \forall k \in \nlcons,\ \ell \leq \primalvector \leq u,\ \primalvector \in \R^{n} \}. \label{eq:nlprelax}
\end{align}
The \emph{mixed integer programming} (MIP) relaxation of~\eqref{eq:minlp} is given by omitting all nonlinear constraints $g_k$ for all $k \in \nlcons$
\begin{align}
  \min \{ c^\T \primalvector\,|\,A\primalvector \geq b,\ \ell \leq \primalvector \leq u,\ \primalvector_i \in \Z\, \forall i \in \intvars \}. \label{eq:miprelax}
\end{align}
Omitting both, integrality requirements and nonlinear constraints, yields the \emph{linear programming} (LP) relaxation of~\eqref{eq:minlp}
\begin{align}
  \min \{ c^\T \primalvector\,|\,A\primalvector \geq b,\ \ell \leq \primalvector \leq u,\ \primalvector \in \R^n\}. \label{eq:lprelax}
\end{align}

All three relaxations provide a lower bound on the optimal solution value of the MINLP~\eqref{eq:minlp}.
MINLP combines discrete decisions and nonlinear functions that are potentially nonconvex.
In theory, linear and convex smooth nonlinear programs are solvable in polynomial time~\cite{khachiyan1979polynomial,vavasis1995complexity}.
In practice, both classes can be solved very efficiently~\cite{bixby2002solving,nocedal2006springer}.
In contrast to that, nonconvexities as imposed by discrete variables or nonconvex nonlinear functions easily lead to
problems that are both \NP-hard in theory and computationally demanding in practice~\cite{vigerske2018scip}.

Commonly used methods to solve convex MINLPs~\eqref{eq:minlp} include
the extended cutting plane algorithm (ECP)~\cite{westerlund1995extended},
the extended supporting hyperplane algorithm~\cite{kronqvist2016extended},
outer approximation (OA)~\cite{duran1986outer,fletcher1994solving}, NLP-based \bandb~\cite{gupta1985branch}, and LP/NLP-based \bandb~\cite{quesada1992lp}.
%For the general case of nonconvex MINLP, the techniques known for convex MINLPs are not applicable in general.
%In that case,
The most commonly used method to solve nonconvex MINLPs is a combination of OA~\cite{kocis1988global,viswanathan1990combined}
and spatial \bandb~\cite{land2010automatic,liberti2003convex,horst2013global}.
Different MINLP solvers either use LP or MIP relaxations or both during the tree search.
For example, \solver{Couenne}~\cite{Couenne} and \scip~\cite{vigerske2018scip} derive valid lower bounds by solving LP relaxations only,
whereas \solver{BARON}~\cite{kilincc2018exploiting,Baron} and \solver{BONMIN}~\cite{bonami2008algorithmic,Bonmin} solve both LP and MIP relaxations.
In contrast to that, only a handful of MINLP solvers provide the possibility to exclusively use NLP relaxations, \eg \solver{BONMIN} and \solver{FICO Xpress Optimizer}~\cite{Xpress}.
For a detailed overview of MINLP solvers that can handle convex and/or nonconvex MINLPs and the implemented algorithm, we refer to~\cite{kronqvist2018review}.\\

In the following, we will focus on MINLP solvers that use a combination of OA and spatial \bandb.
Spatial \bandb is -- analogous to LP-based \bandb~\cite{Dakin1965,LandDoig1960} -- a divide-and-conquer method which splits
the search space sequentially into smaller subproblems that are intended to be easier to solve.
Additionally, convex relaxations are used to compute lower bounds on the individual subproblems.
Based on the computed lower bound, a subproblem can be pruned earlier if the lower bound already exceeds the currently best-known solution.
To divide the search space into smaller pieces, spatial \bandb branches on discrete variables
with a fractional solution value in the relaxation solution.
In addition to that, spatial \bandb uses continuous variables for branching
if they appear in nonconvex terms of nonlinear constraints that are violated by the current relaxation solution.
During this procedure, infeasible subproblems may be encountered.
Infeasibility can either be detected by contradicting variable bounds, derived by domain propagation, or
by an infeasible convex relaxation.
In contrast to modern MIP solvers that can refer to a variety of well-studied techniques,
\eg~\cite{DaveyBolandStuckey2002,SandholmShields2006,achterberg2007constraint}, to 'learn' from infeasible subproblems,
similar techniques for MINLPs exist for certain special cases only.

\section{Conflict Analysis in MINLP}

In this section, we will briefly describe conflict analysis for MIPs of type~\eqref{eq:miprelax}
and the drawbacks when applying these techniques to general MINLP.

\subsection{Technical Background: Conflict Analysis in MIP}
\label{subsec:technicalbackground}
Conflict analysis for MIP has a long history and has its origin in artificial intelligence~\cite{stallman1977forward} and solving satisfiability problems (SAT)~\cite{marques1999grasp}.
Similar ideas are used in constraint programming (CP), see, \eg~\cite{ginsberg1993dynamic,jiang1994no}.
Integrations of these techniques  into MIP were independently
suggested by~\cite{DaveyBolandStuckey2002},~\cite{SandholmShields2006},
and~\cite{achterberg2007constraint}.

If infeasibility is encountered by domain propagation, modern SAT and MIP solvers construct a directed acyclic graph
which represents the logic of how the set of branching decisions led to the detection of infeasibility.
This graph is called the \emph{conflict graph}.
% The vertices of the conflict graph represent bound changes of variables and the arcs $(v,w)$ correspond to bound changes implied by propagation,
% \ie the bound change corresponding to $w$ is based (besides others) on the bound change represented by $v$.
%
% In addition to the inner vertices which represent the bound changes from domain propagation,
% the graph features source vertices for bound changes that correspond to branching decisions and an artificial sink vertex representing the infeasibility.
Valid \emph{conflict constraints} can be derived from cuts in the graph that separate the
branching decisions from an artificial vertex representing the infeasibility.
%Two well-established methods for finding cuts in the conflict graph are the \emph{All-First-Unique-Implication-Point}
%and the \emph{1-First-Unique-Implication-Point} scheme~\cite{zhang2001efficient}.
Based on such a cut, a conflict constraint consists of a set of variables with associated bounds, requiring that in each feasible solution
at least one of the variables has to take a value outside these bounds.

If the LP relaxation of a subproblem with local bounds $\loclb$ and $\locub$ turns out to be infeasible,
it is necessary to identify a set of variables and bound changes that are sufficient to render the infeasibility.
Such a set, the so-called \emph{Farkas proof}~\cite{polik15ismp,WitzigBertholdHeinz2017}, can be constructed by using LP duality theory
that states that exactly one of the systems
\begin{align}
    A\primalvector \geq b,\ \loclb \leq \primalvector \leq \locub & \\
    y^\T A + r^\T \{\loclb,\locub\} = 0,\ y^\T b + r^\T \{\loclb,\locub\} > 0,\ y \geq 0& \label{eq:dualproof}
\end{align}
where $r^\T \{\loclb,\locub\} := \sum_{j \in \allvars\colon r_j > 0} r_j\loclb_j + \sum_{j \in \allvars\colon r_j < 0} r_j\locub_j$,
can be satisfied.
System~\eqref{eq:dualproof} implies a proof of infeasibility \wrt to the local bounds
\begin{align}
    0 < y^\T b + r^\T \{\loclb,\locub\} = y^\T b - (y^\T A)\{\loclb,\locub\} \quad \Longleftrightarrow \quad (y^\T A)\{\loclb,\locub\} < y^\T b.
\end{align}
Consequently, every feasible solution has to satisfy
\begin{align}
    (y^\T A)\primalvector \geq y^\T b, \label{eq:farkasproof}
\end{align}
which is called Farkas proof; it is a globally valid constraint because it is a nonnegative combination of all globally valid constraints.
Thereby, Farkas proofs are a special case of Benders cuts~\cite{Benders1962}.
The Farkas proof is used as a starting point for conflict graph analysis or dual ray analysis.
Note, in MIP conflict graph analysis yields at least one conflict that does not need to be linear, whereas dual ray analysis yields exactly one linear constraint.

\subsection{Conflict Analysis in MINLP}

Only a few publications are dealing with infeasibility in MINLP.
Most of the literature is restricted to a certain class of MINLPs,
\eg conic certificates for convex MINLPs~\cite{coey2018outer} which has been proven to be very successful on
\emph{mixed-integer second-order cone} (MISOCP) problems.
Purely theoretical results for \emph{mixed integer semidefinite programs} (MISDP) were recently published in~\cite{kellner2018irreducible}.
Both publications deal with MINLPs that are infeasible as a whole, and not with the analysis of infeasible subproblems to learn information.
%
% \Todo{ Die connection versteht Timo nicht}
% The probably most general approach for convex MINLPs that is loosely connected to infeasibility analysis in a wider sense comes from the ECP algorithm~\cite{???}.
% The ECP algorithm iteratively improves a polyhedral OA of $g_k(x)$ for all $k \in \nlcons$.
% To strengthen the polyhedral OA the algorithm solves MIP relaxations.
% Every optimal solution $\optMIP{x}$ of the MIP relaxations that violates at least one nonlinear constraint
% is used to separate linear constraints of form
% \begin{align*}
%     g_k(\optMIP{x}) + \nabla g_k(\optMIP{x})\cdot \primalvector \leq 0 \qquad \forall k \in \nlcons\, \colon g_k(\optMIP{x}) > 0,
% \end{align*}
% \ie the first-order Taylor series expansion at $\optMIP{x}$, to stengthen the MIP relaxation.\\

For MINLP algorithms that are based on solving LP relaxations, in
particular, for OA- and ECP-based solvers,
conflict analysis methods for MIP can be applied under certain conditions.
To this end, let us first recap the idea of constructing an LP
relaxation for an MINLP.

%Let $A \in \R^{n \times m}$ be the constraint matrix of all original linear problem constraints
%and assume for all $k \in \nlcons$ that they are nonlinear, i.e., no $g_k$ has a linear representation.
During the tree search, nonlinear functions are approximated by linear functions if they are violated by a relaxation solution.
Let $\tilde{\primalvector}$ be a relaxation solution with $g_k(\tilde{\primalvector}) > 0$.
If $g_k$ is convex, a so-called \emph{gradient cut}
\begin{align*}
    g_k(\tilde{\primalvector}) + \nabla g_k(\tilde{\primalvector})(\primalvector - \tilde{\primalvector}) \leq 0
\end{align*}
is added.
If $g_k$ is nonconvex, convex underestimators are added, see, e.g.,~\cite{vigerske2018scip}.
For quadratic functions, \eg these are the so-called McCormick underestimators~\cite{mccormick1976computability}.
More general nonlinear functions are typically decomposed into functions of a single variable,% (assuming that all $g_k$ are factorable),
for which explicit underestimators are known. Note that gradient cuts are globally valid, while underestimators for non-convex functions
typically involve the local bounds and are hence not globally valid.

For a subproblem $s$ during the tree search, let $\nlconslin^s := \{l^s_1,\ldots,l^s_q\}$ be the index set of all linear approximations of
all $g_k$ with $k \in \nlcons$ that have been added at the node corresponding to $s$ or any of its ancestors.
Hence, it is the current set of (local) \emph{linear relaxation cuts}; all are valid at $s$.
Let $G^s$ be the matrix containing all of these linearizations and $d^s$ be the corresponding right-hand sides.
Thus, the LP relaxation solved for subproblem $s$ reads as
\begin{align}
    \min \{ c^\T \primalvector \;|\; A\primalvector \geq b,\ G^s\primalvector \geq d^s,\ \lb \leq \primalvector \leq \ub\}. \label{eq:lp-oa}
\end{align}
%of the OA is solved at every subproblem $s$ of the \bandb tree.
%Note, the objective function $f(\primalvector)$ is usually replaced by an auxiliary variable $z$ such that $f(\primalvector) \leq z$ holds.
We denote the set of linearizations added at the root node by $\nlconslin^0$.
During the (spatial) \bandb the set of linearizations expands along each path of the tree:
It holds that $\nlconslin^{0} \subseteq \nlconslin^{s_1} \subseteq \ldots \subseteq \nlconslin^{s_p} \subseteq \nlconslin^{s}$ for each path $(0, s_1, \ldots, s_p, s)$.
In analogy to solving MIPs, if~\eqref{eq:lp-oa} is infeasible each ray $(y,w,r)$
in its dual can be used to construct a proof of local infeasibility.
Here, $y_i$ are the dual variables corresponding to $A_{i\cdot}$, $w_l$ are the dual variables corresponding to $G^s_{l \cdot}$ for all $l \in \nlconslin^s$,
and $r_j$ denotes the reduced costs (the duals of the bound constraints) of every variable $\primalvector_j$.
Note that $r_j = c_j - y^\T A_{\cdot j} - w^\T G^s_{\cdot j}$. % denotes the reduced costs of $\primalvector_j$ \wrt~\eqref{eq:lp-oa}.

Hence, a local infeasibility proof \wrt the local bounds $\loclb$ and $\locub$ is given by
\begin{align}
    y^\T b + w^\T d^s + r^\T \{\loclb,\locub\} > 0, \label{eq:oa-dualproof}
\end{align}
In contrast to~\eqref{eq:farkasproof} the constraint $y^\T A\primalvector + w^\T G^s\primalvector \geq y^\T b + w^\T d^s$
is not globally valid in general because linearizations of nonlinear constraints might rely on intermediate local bounds.
Conflict analysis as introduced in~\cite{achterberg2007conflict,WitzigBertholdHeinz2017} only considers globally valid reasons of infeasibility.
Therefore, every local certificate of infeasibility~\eqref{eq:oa-dualproof} needs to be relaxed to consider $\nlconslin^0$ only
\begin{align}
    y^\T b + \bar{w}^\T d^s + \bar{r}^\T \{\loclb,\locub\} > 0, \label{eq:oa-dualproof-relaxed}
\end{align}
where $\bar{w_l} := w_l$, if $l \in \nlconslin^0$, and $\bar{w_l} := 0$, otherwise, and $\bar{r}_j := c_j - y^\T A_{\cdot j} - \bar{w}^\T G^s_{\cdot j}$.
As a consequence, the relaxed certificate~\eqref{eq:oa-dualproof-relaxed}
might not provide an infeasibility proof anymore and cannot be used to generate a conflict constraint.
If, however,~\eqref{eq:oa-dualproof-relaxed} is a valid proof of local infeasibility, all conflict analysis techniques known from MIP can be applied.

\subsection{Locally Valid Certificates of Infeasibility}

% \Todo{HIER MÜSSEN WIR NOCHMAL den Unterschied zwischen conflict graph and dual ray analysis herausarbeiten und was in dem Paper neu ist}
In MIP both conflict graph analysis and dual ray analysis rely on globally valid proofs.
In most MIP solvers, local cuts are applied rarely, if at all. This is very different for non-convex MINLP solvers which rely on local linearization cuts.
%Usually, there is no need to consider locally valid proofs in MIP because local cuts are applied very rarely if at all, this holds at least for \scip\todo{Wie sieht das in Xpress aus? ;-)}.
A computational study within the constraint integer programming and MINLP solver \scip 
%on the impact of MIP techniques on convex and nonconvex
%\emph{mixed-integer quadratically constrained programs} (MIQCPs) showed that disabling conflict graph analysis
%as introduced by~\cite{achterberg2007conflict} would lead to a slowdown of $9\%$~\cite{BertholdGleixnerHeinzetal.2012}.
showed that the impact of conflict graph analysis for general MINLPs is almost negligible~\cite{vigerske2018scip}.
A computational study regarding the impact of dual ray analysis on an MINLP solver
has -- to the best of our knowledge -- never been conducted before.
We present such a computational study in Section~\ref{sec:computational-study}.

The observation that conflict graph analysis on MINLP instances has a much smaller impact than on MIP instances~\cite{BertholdGleixnerHeinzetal.2012,vigerske2018scip}
led to the assumption that a substantial amount of infeasibility proofs of form~\eqref{eq:oa-dualproof-relaxed} were not globally valid.
Hence, they are not suitable for conflict graph analysis as known from the literature and implemented in \scip.
These results indicate that locally added linearization cuts are, non-surprisingly, important to render infeasibility \wrt local bounds.\\

To incorporate local linearizations of nonlinear constraints we propose to generalize dual infeasibility proofs
of subproblem $s$ with local bounds $\loclb$ and $\locub$ as described in Section~\ref{subsec:technicalbackground}
to locally valid certificates of form
\begin{align}
    y^\T b + \hat{w}^\T d^s + \hat{r}^\T \{\loclb,\locub\} > 0, \label{eq:oa-dualproof-local}
\end{align}
incorporating linearizations $\hat{\nlconslin}$ with $\nlconslin^{0} \subseteq \hat{\nlconslin} \subseteq \nlconslin^{s_p}$,
$\hat{w_l} := w_l$, if $l \in \hat{\nlconslin}$, and $\hat{w_l} := 0$, otherwise, and $\hat{r}_j := c_j - y^\T A_{\cdot j} - \hat{w}^\T G^s_{\cdot j}$.
The certificate~\eqref{eq:oa-dualproof-local} is valid for the search tree induced by subproblem $q$,
where $q$ is chosen to satisfy
\begin{align}
q = \min_{q \in \{1,\ldots,s_p\}}\{\nlconslin^{q-1} \subseteq \hat{\nlconslin},\ \hat{\nlconslin} \cap (\nlconslin^{q+1} \setminus \nlconslin^{q}) = \emptyset \}.
\end{align}
Hence, the infeasibility proof might be lifted to an ancestor $q$ of the subproblem $s$ it was created for,
if all local information used for the proof were already available at $q$.
Note that it would be possible to apply conflict graph analysis to~\eqref{eq:oa-dualproof-local}, too.
However, this would introduce a computational overhead because the order of locally apply bound changes
and separated local linearizations needs to be tracked and maintained.
Since conflict graph analysis already comes with an overhead due to maintaining the so-called delta-tree,
\ie complete information about bound deductions and its reasons within the tree,
we omit applying conflict graph analysis on locally valid infeasibility certificates.

\section{Computational Study}
\label{sec:computational-study}

\input{minlp2.tex}

For our computational study, we implemented the generation of locally valid infeasibility certificates
in the academic constraint integer programming solver \scip~\cite{GleixnerBastubbeEifleretal.2018}.
In the following, we refer to \scip with (global) conflict graph analysis as \conflict and \scip with (global) dual ray analysis as \dualray.
Moreover, we refer to \dualray extended by locally valid infeasibility proofs as \dualrayloc.
As a baseline we use \scip with deactivated conflict analysis (\noconflict).
As a test set we use the \MINLP~\cite{MINLPLIB} without instances for which at least one setting finished with numerical violations. %, these instances are excluded.
This yields a test set of 1170 instances.
%The test set consists of 1359 convex and nonconvex MINLPs.
The experiments were run on a cluster of Intel Xeon E5-2690 2.6\,GHz machines with 128\,GB of RAM; a
time limit of $3600$ seconds was set.

Aggregated results of all four settings are shown in Table~\ref{tab:minlp}.
Here,  \bracket{$100$}{tilim} denotes the set of instances for which all settings need at least $100$ seconds
and are solved by at least one setting~\cite{AchterberWunderling2003}.

All settings with activated conflict analysis improve both the running time of \scip, the number of branch-and-bound nodes,
and the number of solved instances.
Moreover, there seems to be a clear ordering: \dualrayloc is superior to \dualray which in turn is superior to \conflict.
%The most amount if instances could be solved by \dualrayloc.
Further, the harder the instances are, the more performance is gained by \dualray and \dualrayloc compared to \conflict.
The number of locally added conflict constraints (\confsloc) by \dualrayloc is on average larger
than the amount of globally added conflict constraints (\confsglb) but in the same order of magnitude.
On the set of nonconvex MINLPs, however, \dualrayloc constructs $11.08$ times more locally than globally valid conflict constraints.
%,which emphasizes the relevance of local linearization cuts.
These results indicate that locally added linearizations of nonlinear constraints are important to render local infeasibility.
%Therefore, ignoring infeasibility certificates incorporating those linearizations leads to drastic loss of information.

When looking into the generation of local proofs into detail,
we could observe that in $5\%$ of all analyzed infeasible LPs no local cut was needed to construct a valid infeasibility certificate, i.e., we could lift
the local conflict to a global one.
%although at least one local linearization cut has a dual solution value different to zero.
For $14\%$ of all local proofs we found a set of local cuts such that $q = \floor{s / 2}$, the conflict could be lifted up at least half of the depth. 
$78\%$ of the local proofs could not be lifted.
Since a lot of infeasibility information is lost, we propose to use a nonlinear relaxation instead.
The theoretical base for nonlinear conflict analysis will be discussed in the following section, whereas the implementation and a computational study is future work.

\section{Outlook and Theoretical Thoughts}

In this final section, we will discuss theoretical considerations how conflict analysis can be directly applied to a nonlinear relaxation of convex MINLPs.
The content described in the following is work-in-progress. % and has not been tested.
%However, it is our motivation for further research in that direction.
At the beginning of this paper, we argued that after LP/MIP-based \bandb, another common method to solve MINLPs is
NLP-based \bandb. We will briefly sketch how a generalization of LP infeasibility analysis can be derived from the KKT-conditions of
convex NLPs.
%solvers use some variant of \bandb, typically based on LP or MIP relaxations.
% In this section, we give some purely theoretical thoughts about relaxations different to LP/MIP and how infeasibility certificates
% can be derived from those relaxations.
Given a convex MINLP of form
% \begin{align}
%     \min_{\primalvector \in X} f(\primalvector) & \notag \\
%     s.t.\quad g_k(\primalvector) &\leq 0 && \forall k \in \nlcons, \label{eq:convex-minlp}\\
%     h_e(\primalvector) &= 0 && \forall e \in \nlconseq, \notag
% \end{align}
\begin{align}
    \min_{\primalvector \in X} \{ f(\primalvector) \;|\; g_k(\primalvector) \leq 0\; \forall k \in \nlcons,\ h_e(\primalvector) = 0\ \forall e \in \nlconseq\}, \label{eq:convex-minlp}
\end{align}
where $f, g_k$ are convex, continuously differentiable functions over $\R^n$ and $h_e$ are affine functions.
For every optimal solution $\optMIP{\primalvector}$ of~\eqref{eq:convex-minlp} of the
(convex) NLP relaxation of~\eqref{eq:convex-minlp} there exist $\lambda \geq 0$ such that it holds that
\begin{align}
    \nabla f(\optMIP{\primalvector}) + \sum_{k \in \nlcons} \lambda_{k}\nabla g_k(\optMIP{\primalvector}) + \sum_{e \in \nlconseq} \mu_{e}\nabla h_e(\optMIP{\primalvector}) = 0 \label{eq:derivative-lagrangian}, \quad
    \lambda_k g_k(\optMIP{\primalvector}) = 0.
\end{align}
These conditions raise from the so-called \emph{Karush-Kuhn-Tucker-Conditions}~\cite{kuhn2014nonlinear}.
Equality~\eqref{eq:derivative-lagrangian} is the gradient of the \emph{Lagrangian dual} that reads as
\begin{align}
    \Lagrangian(\primalvector,\lambda,\mu) := f(\primalvector) + \sum_{k \in \nlcons} \lambda_k g_k(\primalvector) + \sum_{e \in \nlconseq} \mu_e h_e(\primalvector),
\end{align}
with $\lambda \geq 0$ and $\mu \in \R^{\Abs{\nlconseq}}$.
By duality theory, the \emph{Lagrangian dual function} which reads as
$\LagrangianFunktion(\lambda, \mu) := \sup_{\lambda,\mu} \Lagrangian(\primalvector, \lambda, \mu)$
yields a lower bound on the optimal value of~\eqref{eq:convex-minlp}. Maximizing $\LagrangianFunktion(\lambda, \mu)$ would give
the tightest lower bound of~\eqref{eq:convex-minlp}, and strict duality of convex optimization tells us that this is equivalent to the optimal
value of~\eqref{eq:convex-minlp}. %If the (convexNLP~\eqref{eq:convex-minlp}
Consequently, if there exists $(\lambda^{\star}, \mu^{\star})$ such that
% \begin{align}
%     \sum_{k \in \nlcons} \lambda^{\star}_k g_k(\primalvector) + \sum_{e \in \nlconseq} \mu^{\star}_{e} h_e(\primalvector) > 0,
% \end{align}
$\sum_{k \in \nlcons} \lambda^{\star}_k g_k(\primalvector) + \sum_{e \in \nlconseq} \mu^{\star}_{e} h_e(\primalvector) > 0$,
then the dual is unbounded and thus $(\lambda^{\star}, \mu^{\star})$ proofs infeasibility of~\eqref{eq:convex-minlp}.
Even though Slater regularity does not hold for infeasible points
\footnote{If one wanted to assume
regularity on the constraint functions of~\eqref{eq:convex-minlp}, linear independence constraint classification would be applicable.},
\begin{align}
    \sum_{k \in \nlcons} \lambda^{\star}_k g_k(\primalvector) + \sum_{e \in \nlconseq} \mu^{\star}_{e} h_e(\primalvector) \leq 0 \label{eq:nlp-farkas}
\end{align}
is a valid inequality for~\eqref{eq:convex-minlp}; it is a convex combination (defined by the dual multipliers)
of the constraints of~\eqref{eq:convex-minlp}.
Inequality~\eqref{eq:nlp-farkas} is the convex optimization equivalent of the Farkas proof~\eqref{eq:farkasproof}.

Assume that constraint~\eqref{eq:nlp-farkas} is given as proof of infeasibility for a subproblem within an NLP-based \bandb.
%\wrt to the current local bounds.
If no local cuts are involved in the infeasibility proof,
inequality~\eqref{eq:nlp-farkas} is a globally valid convex nonlinear constraint.
Note in this context that gradient cuts are globally valid.
%The duality theory how a certificate of in to otained is very close the one used in MIP.
%Principally, the local bounds involved in the infeasibility proof could be used to seed a conflict graph analysis.
%Like in the linear case, the problem is that the initial reason might 

Clearly, inequality~\eqref{eq:nlp-farkas} holds for all non-negative
$\lambda^{\star}$. The following observation makes
the concrete $(\lambda^{\star}, \mu^{\star})$ from the infeasibility
proof interesting to use as global information inside a \bandb tree
search for convex MINLP. Consider the linearization at an infeasible point $\primalvector^{\star}$
\begin{align}
  \nabla g_k(\primalvector^\star)^\T  (\primalvector - \primalvector^\star) \leq 0 \quad \Leftrightarrow \quad  \nabla g_k(\primalvector^\star)^\T  \primalvector \leq \nabla g_k(\primalvector^\star)^\T  \primalvector^\star\; \forall k \in \nlcons.
\end{align}
Then, the corresponding dual multipliers $\lambda^{\star}$ give the (linear)
Farkas proof
\begin{align}
\sum_{k \in \nlcons} \lambda^{\star}_k \nabla g_k(\primalvector^\star) + \sum_{e \in \nlconseq} \mu^{\star}_e \nabla h_e(\primalvector^\star) &= 0 \\
\sum_{k \in \nlcons} \lambda^{\star}_k \nabla g_k(\primalvector^\star)^\T \primalvector^\star + \sum_{e \in \nlconseq} \mu^{\star}_e \nabla h_e(\primalvector^\star)^\T  \primalvector^\star &< 0.
\end{align}
Hence, as in the case of dual ray analysis for MIP,
inequality~\eqref{eq:nlp-farkas} is a single inequality that would
have provided the infeasibility proof from its derivative. The hope
(which is true for MIP) is that it is a good candidate to detect
infeasibility by propagation (under the use of integrality
information) in other parts of the search tree, and might be a
meaningful aggregation of problem constraints to create cuts from.

For many NLP solvers, in particular \emph{dual active set
  methods}~\cite{murty1988linear,nocedal2006nonlinear,forsgren2016primal}
and barrier algorithms~\cite{Mehrotra1992,Meszaros1999,Wachter2009}, dual multipliers will be readily available.
The added advantage of active set methods is that they typically yield a sparse dual weight-vector $(\lambda,\mu)$.
This might come in handy when the local bounds involved in the infeasibility proof should be used to seed a conflict graph analysis.
Like in the linear case, the problem is that the initial reason will typically be too large to be meaningful.

All of this is subject to further investigation.
We plan to implement NLP-based conflict analysis into the academic constraint integer programming solver \scip
and to study its impact on solver behavior.
As in the MIP case, infeasibility information might be used in several other contexts, consider hybrid branching~\cite{AchterbergBerthold2009},
conflict-driven diving heuristics~\cite{WitzigGleixner2018}, and
also rapid learning~\cite{BertholdFeydyStuckey2010,BertholdStuckeyWitzig2018}.
% \todo{rausfinden, was die klassischen Zitate für KKT und Active Set sind}

\subsection*{Acknowledgments}

We thank Zsolt Csizmadia for his valuable comments on Section 4. The work for this article has been conducted within the Research Campus Modal
funded by the German Federal Ministry of Education and Research (fund number 05M14ZAM).
We thank three anonymous reviewers for their valuable suggestions and helpful comments.

\bibliographystyle{abbrv}
\bibliography{Bibliography}

% \ifthenelse{\zibreport=1}
% {
% \clearpage
% \pagebreak
%
% \begin{appendix}
%
% \section{Appendix}
%
%
% \end{appendix}
% }
% {}

\end{document}

%% file: Defines.tex
\usepackage{amsmath}
\usepackage{amssymb}

\usepackage{xspace}
\usepackage{xcolor}

\usepackage{ifthen}
\usepackage{xifthen}

\usepackage[colorlinks=false, pdftex]{hyperref}
\hypersetup{
    colorlinks,
    citecolor={blue!80!black},
    linkcolor={blue!80!black},
    urlcolor={blue!80!black}
}

\usepackage[textsize=small,textwidth=3.3cm]{todonotes}

\newcolumntype{L}{>{\raggedright\arraybackslash}X}%
\newcolumntype{C}{>{\centering\arraybackslash}X}%
\newcolumntype{R}{>{\raggedleft\arraybackslash}X}%

\newcommand{\cleaninst}{all}

\newcommand{\bracket}[2]{[#1,#2]}

\newcommand{\solved}{\textbf{solved}\xspace}
\newcommand{\nodes}{\textbf{nodes}\xspace}
\newcommand{\nodesQ}{\textbf{nodes\textsubscript{Q}}\xspace}
\renewcommand{\time}{\textbf{time}\xspace}
\newcommand{\timeQ}{\textbf{time\textsubscript{Q}}\xspace}
\newcommand{\confsloc}{\textbf{confs\textsubscript{loc}}\xspace}
\newcommand{\confsglb}{\textbf{confs\textsubscript{glb}}\xspace}

\newcommand{\MINLP}{\textsc{MINLPLIB}\xspace}

\newcommand{\setting}[1]{\texttt{#1}}
\newcommand{\noconflict}{\setting{noconflict}\xspace}
\newcommand{\conflict}{\setting{confgraph}\xspace}
\newcommand{\dualray}{\setting{dualray}\xspace}
\newcommand{\dualrayloc}{\setting{dualray-loc}\xspace}

\newcommand{\solver}[1]{\texttt{#1}}
\newcommand{\scip}{\solver{SCIP}\xspace}

\newcommand{\ie}{i.e.,\xspace}
\newcommand{\eg}{e.g.,\xspace}

\newcommand{\wrt}{w.r.t.\xspace}

\newcommand{\bandb}{branch-and-bound\xspace}

\newcommand{\lb}{\ell}
\newcommand{\ub}{u}
\newcommand{\loclb}{\lb^\prime}
\newcommand{\locub}{\ub^\prime}

\newcommand{\NP}{$\mathcal{NP}$}

\newcommand{\Abs}[1]{\vert #1\vert}

\newcommand{\intvars}{\mathcal{I}}
\newcommand{\allvars}{\mathcal{N}}

\newcommand{\optMIP}[1]{#1^{\star}}

\newcommand{\floor}[1]{\lfloor #1\rfloor}

\newcommand{\minact}[4]{\ifthenelse{\isempty{#4}}{\Delta_{\min}(#1,#2,#3)}{\Delta_{\min}^{#4}(#1,#2,#3)}}
\newcommand{\maxact}[4]{\ifthenelse{\isempty{#4}}{\Delta_{\max}(#1,#2,#3)}{\Delta_{\max}^{#4}(#1,#2,#3)}}

\newcommand{\Z}{\mathbb{Z}\xspace}
\newcommand{\R}{\mathbb{R}\xspace}

\newcommand{\nlcons}{\mathcal{K}\xspace}
\newcommand{\nlconseq}{\mathcal{E}\xspace}
\newcommand{\nlconslin}{\mathcal{G}\xspace}

\newcommand{\Lagrangian}{\mathcal{L}\xspace}
\newcommand{\LagrangianFunktion}{q\xspace}

\newcommand{\var}{\ensuremath{x}\xspace}
\newcommand{\primalvector}{\var}
\newcommand{\T}{\mathsf{T}}

%% file: minlp2.tex
\begin{table}[t]
\scriptsize
\caption{Aggregated results on \MINLP}
\label{tab:minlp}
\begin{tabularx}{\textwidth}{lRRRRRRrR}
\toprule
                      &    \#  &  \solved  &    \time &    \nodes &   \timeQ  & \nodesQ  &   \confsglb &  \confsloc \\
\midrule
\cleaninst            &        &           &          &           &           &          &             &            \\
\noconflict           &  1170  &       689 &    79.11 &   3014.25 &    1.000  &   1.000  &          -- &       --   \\
\conflict             &  1170  &       694 &    77.94 &   2952.07 &    0.985  &   0.979  &     9679.01 &       --   \\
\dualray              &  1170  &       695 &    76.78 &   2871.86 &    0.970  &   0.953  &     1359.92 &       --   \\
\dualrayloc           &  1170  &   \bf{698}&    76.35 &   2841.90 &\bf{0.965} &\bf{0.943}&     1338.65 &  3192.50   \\[.25em]
\bracket{100}{tilim}  &        &           &          &           &           &          &             &            \\
\noconflict           &    99  &        83 &   638.34 &  86860.54 &    1.000  &   1.000  &          -- &       --   \\
\conflict             &    99  &        88 &   563.06 &  74251.69 &    0.882  &   0.855  &    23653.88 &       --   \\
\dualray              &    99  &        89 &   458.28 &  62890.08 &    0.718  &   0.724  &     2019.46 &       --   \\
\dualrayloc           &    99  &    \bf{92}&   429.31 &  59629.05 &\bf{0.673} &\bf{0.686}&     2086.62 &  3177.67   \\%[.25em]
% \bracket{1000}{tilim} &        &           &          &           &           &          &             &            \\
% \noconflict           &    45  &        29 &  1841.77 & 209594.28 &    1.000  &   1.000  &          -- &       --   \\
% \conflict             &    45  &        34 &  1591.92 & 175862.89 &    0.864  &   0.839  &    42899.80 &       --   \\
% \dualray              &    45  &        35 &   998.40 & 114341.05 &    0.542  &   0.546  &     2012.24 &       --   \\
% \dualrayloc           &    45  &        38 &   875.54 & 103595.17 &    0.475  &   0.494  &     2038.60 &  5265.95   \\
\bottomrule
\end{tabularx}
\end{table}

%% file: main.bbl
\begin{thebibliography}{10}

\bibitem{achterberg2007conflict}
T.~Achterberg.
\newblock Conflict analysis in mixed integer programming.
\newblock {\em Discrete Optimization}, 4(1):4--20, 2007.

\bibitem{achterberg2007constraint}
T.~Achterberg.
\newblock Constraint integer programming, 2007.

\bibitem{AchterbergBerthold2009}
T.~Achterberg and T.~Berthold.
\newblock Hybrid branching.
\newblock In W.-J. van Hoeve and J.~N. Hooker, editors, {\em Integration of AI
  and OR Techniques in Constraint Programming for Combinatorial Optimization
  Problems, 6th International Conference, CPAIOR 2009}, volume 5547 of {\em
  Lecture Notes in Computer Science}, pages 309--311. Springer Berlin
  Heidelberg, May 2009.

\bibitem{AchterberWunderling2003}
T.~Achterberg and R.~Wunderling.
\newblock Mixed integer programming: Analyzing 12 years of progress.
\newblock In M.~J{\"u}nger and G.~Reinelt, editors, {\em Facets of
  Combinatorial Optimization: Festschrift for Martin Gr{\"o}tschel}, pages
  449--481. Springer Berlin Heidelberg, 2013.

\bibitem{Baron}
{BARON}.
\newblock \url{https://minlp.com/baron}.

\bibitem{Benders1962}
J.~F. Benders.
\newblock Partitioning procedures for solving mixed-variables programming
  problems.
\newblock {\em {Numerische Mathematik}}, 4(1):238--252, 1962.

\bibitem{BertholdFeydyStuckey2010}
T.~Berthold, T.~Feydy, and P.~J. Stuckey.
\newblock Rapid learning for binary programs.
\newblock In A.~Lodi, M.~Milano, and P.~Toth, editors, {\em Integration of AI
  and OR Techniques in Constraint Programming for Combinatorial Optimization
  Problems, 7th International Conference, CPAIOR 2010}, volume 6140 of {\em
  LNCS}, pages 51--55. Springer Berlin Heidelberg, June 2010.

\bibitem{BertholdGleixnerHeinzetal.2012}
T.~Berthold, A.~M. Gleixner, S.~Heinz, and S.~Vigerske.
\newblock Analyzing the computational impact of {MIQCP} solver components.
\newblock {\em Numerical Algebra, Control and Optimization}, 2(4):739 -- 748,
  2012.

\bibitem{BertholdStuckeyWitzig2018}
T.~Berthold, P.~J. Stuckey, and J.~Witzig.
\newblock Local rapid learning for integer programs.
\newblock Technical Report 18-56, ZIB, Takustr.~7, 14195 Berlin, 2018.

\bibitem{bixby2002solving}
R.~E. Bixby.
\newblock Solving real-world linear programs: A decade and more of progress.
\newblock {\em Operations research}, 50(1):3--15, 2002.

\bibitem{bonami2008algorithmic}
P.~Bonami, L.~T. Biegler, A.~R. Conn, G.~Cornu{\'e}jols, I.~E. Grossmann, C.~D.
  Laird, J.~Lee, A.~Lodi, F.~Margot, N.~Sawaya, et~al.
\newblock An algorithmic framework for convex mixed integer nonlinear programs.
\newblock {\em Discrete Optimization}, 5(2):186--204, 2008.

\bibitem{Bonmin}
{Bonmin}.
\newblock \url{https://projects.coin-or.org/Bonmin}.

\bibitem{coey2018outer}
C.~Coey, M.~Lubin, and J.~P. Vielma.
\newblock Outer approximation with conic certificates for mixed-integer convex
  problems.
\newblock {\em arXiv preprint arXiv:1808.05290}, 2018.

\bibitem{Couenne}
{Couenne}.
\newblock \url{https://www.coin-or.org/Couenne/}.

\bibitem{Dakin1965}
R.~J. Dakin.
\newblock A tree-search algorithm for mixed integer programming problems.
\newblock {\em The Computer Journal}, 8(3):250--255, 1965.

\bibitem{DaveyBolandStuckey2002}
B.~Davey, N.~Boland, and P.~J. Stuckey.
\newblock Efficient intelligent backtracking using linear programming.
\newblock {\em INFORMS Journal of Computing}, 14(4):373--386, 2002.

\bibitem{duran1986outer}
M.~A. Duran and I.~E. Grossmann.
\newblock An outer-approximation algorithm for a class of mixed-integer
  nonlinear programs.
\newblock {\em Mathematical programming}, 36(3):307--339, 1986.

\bibitem{Xpress}
{FICO Xpress Optimizer}.
\newblock \url{https://www.fico.com/de/products/fico-xpress-optimization}.

\bibitem{fletcher1994solving}
R.~Fletcher and S.~Leyffer.
\newblock Solving mixed integer nonlinear programs by outer approximation.
\newblock {\em Mathematical programming}, 66(1-3):327--349, 1994.

\bibitem{forsgren2016primal}
A.~Forsgren, P.~E. Gill, and E.~Wong.
\newblock Primal and dual active-set methods for convex quadratic programming.
\newblock {\em Mathematical programming}, 159(1-2):469--508, 2016.

\bibitem{ginsberg1993dynamic}
M.~L. Ginsberg.
\newblock Dynamic backtracking.
\newblock {\em Journal of Artificial Intelligence Research}, 1:25--46, 1993.

\bibitem{GleixnerBastubbeEifleretal.2018}
A.~Gleixner, M.~Bastubbe, L.~Eifler, T.~Gally, G.~Gamrath, R.~L. Gottwald,
  G.~Hendel, C.~Hojny, T.~Koch, M.~E. L{\"u}bbecke, S.~J. Maher,
  M.~Miltenberger, B.~M{\"u}ller, M.~E. Pfetsch, C.~Puchert, D.~Rehfeldt,
  F.~Schl{\"o}sser, C.~Schubert, F.~Serrano, Y.~Shinano, J.~M. Viernickel,
  M.~Walter, F.~Wegscheider, J.~T. Witt, and J.~Witzig.
\newblock {The SCIP Optimization Suite 6.0}.
\newblock Technical Report 18-26, ZIB, Takustr. 7, 14195 Berlin, 2018.

\bibitem{gupta1985branch}
O.~K. Gupta and A.~Ravindran.
\newblock Branch and bound experiments in convex nonlinear integer programming.
\newblock {\em Management science}, 31(12):1533--1546, 1985.

\bibitem{horst2013global}
R.~Horst and H.~Tuy.
\newblock {\em Global optimization: Deterministic approaches}.
\newblock Springer Science \& Business Media, 2013.

\bibitem{jiang1994no}
Y.~Jiang, T.~Richards, and B.~Richards.
\newblock No-good backmarking with min-conflict repair in constraint
  satisfaction and optimization.
\newblock In {\em PPCP}, volume~94, pages 2--4. Citeseer, 1994.

\bibitem{kellner2018irreducible}
K.~Kellner, M.~E. Pfetsch, and T.~Theobald.
\newblock Irreducible infeasible subsystems of semidefinite systems.
\newblock {\em arXiv preprint arXiv:1804.01327}, 2018.

\bibitem{khachiyan1979polynomial}
L.~G. Khachiyan.
\newblock A polynomial algorithm in linear programming.
\newblock In {\em Doklady Academii Nauk SSSR}, volume 244, pages 1093--1096,
  1979.

\bibitem{kilincc2018exploiting}
M.~R. K{\i}l{\i}n{\c{c}} and N.~V. Sahinidis.
\newblock {Exploiting integrality in the global optimization of mixed-integer
  nonlinear programming problems with BARON}.
\newblock {\em Optimization Methods and Software}, 33(3):540--562, 2018.

\bibitem{kocis1988global}
G.~R. Kocis and I.~E. Grossmann.
\newblock {Global optimization of nonconvex mixed-integer nonlinear programming
  (MINLP) problems in process synthesis}.
\newblock {\em Industrial \& engineering chemistry research}, 27(8):1407--1421,
  1988.

\bibitem{kronqvist2018review}
J.~Kronqvist, D.~Bernal, A.~Lundell, and I.~Grossmann.
\newblock A review and comparison of solvers for convex {MINLP}.
\newblock {\em Preprint, Optimization Online}, 2018.

\bibitem{kronqvist2016extended}
J.~Kronqvist, A.~Lundell, and T.~Westerlund.
\newblock The extended supporting hyperplane algorithm for convex mixed-integer
  nonlinear programming.
\newblock {\em Journal of Global Optimization}, 64(2):249--272, 2016.

\bibitem{kuhn2014nonlinear}
H.~W. Kuhn and A.~W. Tucker.
\newblock Nonlinear programming.
\newblock In {\em Traces and emergence of nonlinear programming}, pages
  247--258. Springer, 2014.

\bibitem{LandDoig1960}
A.~H. Land and A.~G. Doig.
\newblock An automatic method of solving discrete programming problems.
\newblock {\em Econometrica}, 28(3):497--520, 1960.

\bibitem{land2010automatic}
A.~H. Land and A.~G. Doig.
\newblock An automatic method for solving discrete programming problems.
\newblock In {\em 50 Years of Integer Programming 1958-2008}, pages 105--132.
  Springer, 2010.

\bibitem{liberti2003convex}
L.~Liberti and C.~C. Pantelides.
\newblock Convex envelopes of monomials of odd degree.
\newblock {\em Journal of Global Optimization}, 25(2):157--168, 2003.

\bibitem{marques1999grasp}
J.~P. Marques-Silva and K.~Sakallah.
\newblock Grasp: A search algorithm for propositional satisfiability.
\newblock {\em Computers, IEEE Transactions on}, 48(5):506--521, 1999.

\bibitem{mccormick1976computability}
G.~P. McCormick.
\newblock {Computability of global solutions to factorable nonconvex programs:
  Part I—Convex underestimating problems}.
\newblock {\em Mathematical programming}, 10(1):147--175, 1976.

\bibitem{Mehrotra1992}
S.~Mehrotra.
\newblock On the implementation of a primal-dual interior point method.
\newblock {\em SIAM Journal on optimization}, 2(4):575--601, 1992.

\bibitem{Meszaros1999}
C.~M{\'e}sz{\'a}ros.
\newblock The {BPMPD} interior point solver for convex quadratic problems.
\newblock {\em Optimization Methods and Software}, 11(1-4):431--449, 1999.

\bibitem{MINLPLIB}
{MINLPLib}.
\newblock {Githash 033934c0}.
\newblock \url{http://www.minlplib.org/}.

\bibitem{murty1988linear}
K.~G. Murty and F.-T. Yu.
\newblock {\em Linear complementarity, linear and nonlinear programming},
  volume~3.
\newblock Citeseer, 1988.

\bibitem{nocedal2006springer}
J.~Nocedal and S.~Wright.
\newblock Springer series in operations research and financial engineering.
\newblock {\em Numerical Optimization, 2nd edn. Springer, Berlin}, 2006.

\bibitem{nocedal2006nonlinear}
J.~Nocedal and S.~J. Wright.
\newblock {\em Nonlinear Equations}.
\newblock Springer, 2006.

\bibitem{polik15ismp}
I.~P\'olik.
\newblock {Some more ways to use dual information in MILP}.
\newblock In {\em International Symposium on Mathematical Programming},
  Pittsburgh, PA, 2015.

\bibitem{quesada1992lp}
I.~Quesada and I.~E. Grossmann.
\newblock An {LP/NLP} based branch and bound algorithm for convex {MINLP}
  optimization problems.
\newblock {\em Computers \& chemical engineering}, 16(10-11):937--947, 1992.

\bibitem{SandholmShields2006}
T.~Sandholm and R.~Shields.
\newblock Nogood learning for mixed integer programming.
\newblock In {\em Workshop on Hybrid Methods and Branching Rules in
  Combinatorial Optimization, Montr{\'e}al}, 2006.

\bibitem{stallman1977forward}
R.~M. Stallman and G.~J. Sussman.
\newblock Forward reasoning and dependency-directed backtracking in a system
  for computer-aided circuit analysis.
\newblock {\em Artificial intelligence}, 9(2):135--196, 1977.

\bibitem{vavasis1995complexity}
S.~A. Vavasis.
\newblock Complexity issues in global optimization: a survey.
\newblock In {\em Handbook of global optimization}, pages 27--41. Springer,
  1995.

\bibitem{vigerske2018scip}
S.~Vigerske and A.~Gleixner.
\newblock {SCIP: Global optimization of mixed-integer nonlinear programs in a
  branch-and-cut framework}.
\newblock {\em Optimization Methods and Software}, 33(3):563--593, 2018.

\bibitem{viswanathan1990combined}
J.~Viswanathan and I.~E. Grossmann.
\newblock A combined penalty function and outer-approximation method for
  {MINLP} optimization.
\newblock {\em Computers \& Chemical Engineering}, 14(7):769--782, 1990.

\bibitem{Wachter2009}
A.~W{\"a}chter.
\newblock Short tutorial: getting started with {Ipopt} in 90 minutes.
\newblock In {\em Dagstuhl Seminar Proceedings}. Schloss
  Dagstuhl-Leibniz-Zentrum f{\"u}r Informatik, 2009.

\bibitem{westerlund1995extended}
T.~Westerlund and F.~Pettersson.
\newblock An extended cutting plane method for solving convex {MINLP} problems.
\newblock {\em Computers \& Chemical Engineering}, 19:131--136, 1995.

\bibitem{WitzigBertholdHeinz2017}
J.~Witzig, T.~Berthold, and S.~Heinz.
\newblock Experiments with conflict analysis in mixed integer programming.
\newblock In {\em Integration of AI and OR Techniques in Constraint Programming
  for Combinatorial Optimization Problems, 14th International Conference,
  CPAIOR 2017}, volume 10335 of {\em LNCS}, pages 211--220. Springer Berlin
  Heidelberg, May 2017.

\bibitem{WitzigGleixner2018}
J.~Witzig and A.~Gleixner.
\newblock {Conflict-Driven Heuristics for Mixed Integer Programming}.
\newblock In preparation.

\end{thebibliography}
